\documentclass{amsart}
\usepackage{geometry}                
\usepackage{graphicx}
\usepackage{amssymb}
\usepackage{epstopdf,amsmath}
\usepackage{enumitem}

\newtheorem*{lemma*}{Lemma}

\newtheorem*{prop*}{Proposition}
\newtheorem*{theorem*}{Theorem}

\title{What math means to me}

\author{Larry Guth}

\begin{document}

\maketitle


Studying math is a big part of my life.   This essay is about how I see math, what I value about it,  and what it means to me.

I had a couple motivations to write it.   One motivation is related to AI.   Developments in AI are prompting the math community to think about what we care about and what directions we would like to go in the future.   
At the ICM this summer,  several speakers encouraged all of us to reflect on our goals and values.   But I've also thought for a long time about trying to write about my personal experiences doing math.    The current situation gave me the extra motivation to do it.


\section{Things come together}

When I talk with people who aren't mathematicians, they often describe their math experiece as an easy part, when they were in grade school,  and then a confusing incomprehensible part starting at some later time.  The ``easy part'' includes counting,  adding, subtracting, ...  Different people experience a transition from easy to incomprehensible at different spots,  and they usually describe the transition as abrupt.

I don't see math as divided into parts in this way.  To begin discussing how I see math, I think that counting is a great example.   Learning to count is not easy -- I think we just don't remember what it was like.   

I remember when my son was learning to count.   He was able to recite numbers up to 30 or so.   We were sitting with a tape measure,  and he wanted to understand which number goes with which line on the tape measure,  and it was very difficult for him.   

There is a time when a child is able to recite the numbers from one to ten,  but they can't figure out how many spoons to put on the dinner table so that each person will have one.  A lot has to come together for a child to count spoons for the dinner table.   The group of people having dinner, and the group of spoons in the drawer,  and the numbers we recite are all related,  and the child has to understand all those relationships.   


Learning to count comes from a special period of childhood when the world keeps opening up as the child understands more.  Most of us remember this time only dimly.   A special thing about studying math is that we have experiences like that all our lives.

When we learn math,  things come together.  

Pictures and numbers and equations come together when we learn equations for lines and circles.   I think this connection is really remarkable.   Lines and circles are fundamental shapes,  used by people for thousands of years,  drawn by children in kindergarten.   Numbers and equations are fundamental too,  but they seem from a different world.   Children learn multiplication to answer questions like,  ``if there are eight children in the class and we give them each four crackers,  how many crackers do we need to take out?''   I don't know how the idea of multiplication first developed historically,  but perhaps for similar reasons.   It is surprising that this cracker counting problem is related to lines and even more surprising that it is related to circles.   And yet circles and lines are described by elegant equations based on multiplication and addition.   This connection lets us ``look at'' circles and lines in a completely different way.  A lot has to come together to understand the equation for a circle.  

In calculus,  we connect the formula for a graph $y = f(x)$ and the way the graph looks.   In particular,  we learn how the slope of the graph is related to the formula for the function $f$.   Later we learn how the area under a graph is related to  the idea of slope.   A lot has to come together for a student to understand how to use calculus to find the area under a parabola.

Many areas of advanced mathematics are named for two things coming together: algebraic topology,  differential geometry,  algebraic geometry...  The first branch of math I studied as a PhD student was called quantitative topology -- meaning quantitative estimates and topology coming together.   

The next branch of math I studied was Fourier analysis.   In Fourier analysis,  we learn that any function can be represented in two different ways,  which are called the physical space representation and the frequency space representation.    Fourier analysis is about looking at functions from these two very different perspectives.   Things that are complicated from one perspective become clearer from the other.   By combining the two perspectives,  we are able to understand things better.   My graduate students and I have spent many hours together at the board drawing pictures of functions in frequency space and trying to understand what they look like in physical space -- trying to see how these two points of view fit together.   

Fourier analysis is helpful for understanding things that oscillate like a pendulum or a vibrating string.   It was originally developed to help understand physical systems such as the behavior of heat or waves.   But it has turned out to help understanding other very different things like the distribution of prime numbers or the structure of crystals or certain fast algorithms on classical and quantum computers.   A lot has to come together for a mathematician to understand how vibrating strings relate to crystals,  prime numbers,  and fast algorithms.   In many other areas,  mathematicians can tell similar stories.

As we spend time getting to know questions,  we see more and more relationships between them.   These relationships can help to answer questions.   And good questions can wake us up to see relationships we hadn't seen before.    The mathematical world is more connected than I expected as a young student.  As we keep studying,  things fit together into a larger whole.  This connection and wholeness is a remarkable feature of the mathematical world -- maybe of nature.   When we practice math,  we cultivate it.

\section{Confusion and frustration and doubt and being wrong}

In the work we do day to day,  we are only rarely aware of that sense of connection and wholeness.   Our immediate experience is usually full of confusion,  frustration,  and doubt. 

I have crumpled papers,  scribbled over text,  and jammed my pencil deep into pads.   I have gotten excited about some idea and then realized it didn't make sense because of something I did wrong three weeks ago.   I remember as a graduate student talking with my parents on the phone and telling them I felt like Wiley E. Coyote.  My experience of doing research often reminds me of Wiley E.  Coyote: cooking up an elaborate strategy,  setting it in motion,  watching it go awry,  and then getting smacked in the head by a bowling ball or discovering I have walked off a cliff.  

Research can feel like a roller coaster.  While I'm excited about an idea, I sometimes fantasize about how remarkable people will think the work is -- especially if I'm working on a kind of famous problem like the Kakeya problem.  When I realize the idea doesn't make sense, I start to have doubts about whether I'm ``good enough'' or about whether the work matters.  

These experiences are a big part of practicing math.   And out of all these difficult experiences,  I think the experience of being wrong is especially important.   We try like crazy to do something,  and we come up with an idea,  and we really want it to be true,  and we realize that it's wrong,  and we admit it.

This hard work and hard experience can teach us to tell sense from non-sense.  
It can teach us that nature is how it is and that it doesn't care what we wish was true.    It can teach us to be honest with ourselves.   

The experience of being wrong can come in all kinds of work and study,  but the experience works a little differently in different areas.   In math,  we have worked hard to make true and false particularly clear and sharply defined,  and I think that helps force us to face being wrong.   

Doing math,  I've faced being wrong in a way that I didn't always make myself face it in other things I studied.   For instance,  I think that my experience doing math has made me a somewhat better writer.   When I was in high school and college,  I really liked writing.   I would try to write a sentence that was interesting and out-of-the-box.   Then I would try to add another interesting sentence.   These days,  after a long time writing about math,  I write the first sentence and I look at it.  ``Do I really think this is true?''  Usually not.   Cross it out.  Start again.  



\section{Asking questions}

Asking questions is a big part of doing math.  
Questions can point out our confusions.   They can highlight an example that doesn't fit our intuition.  They can wake us up to new connections.   

When we get stuck on a question, we try to ask a more basic question.  And when we get stuck on that basic question,  we try to ask even more basic questions.   I once heard a saying about this process that went like this: ``if there is something I don't understand,  there is usually something simpler that I also don't understand.''  I am not sure exactly of the history of this saying,  but it comes to my mind often.   For me, this philosophy is an important part of mathematics.

This constant attempt to ask more basic questions can help us to see connections between things.   We keep an eye out for anything that seems arbitrary.   If some feature of the question seems arbitrary,  we try to generalize that feature and ask a more basic, more general question.   Question A and Question B may each have all sorts of arbitrary aspects which look very different.   But if we thoughtfully remove the arbitrary aspects and look underneath,  then the questions may start to look the same.   

I used to teach a one day undergrad seminar about problem solving called ``Towards a simpler problem''.   The point of the seminar was that,  faced with a really hard problem,  when we are totally stuck,  we can turn to a simpler problem that sheds some light on the original question.   I wanted to show students that,  faced with,  say,  the problem in Fermat's last theorem,  they don't have to be frozen and overwhelmed,  but that they can start to engage with it.   They can ask smaller questions and start to figure those out and that helps understand the bigger question better and that leads to more questions....  

Questions are important because when we are faced with a problem too big for us to budge,  we can still engage with it by asking questions.   

Faced with a hard problem,  one approach is to look for a problem that is halfway between what we already understand and our goal.   But finding such a `halfway problem' is often not easy.   There is a well known issue in biology that the simplest organism is a single cell,  which is extremely complicated.   There is no such thing as ``half a cell'' to study as a warmup problem.   Many problems in math feel somewhat like this -- it's not at all clear how to divide the problem in half.   Nevertheless,  in creative ways,  mathematicians find simpler problems that shed light on the original problem.

Here is an example from my personal experience.   The restriction problem is a difficult problem in Fourier analysis or PDE about the interference patterns that appear when we superimpose many waves with different frequencies.   Mathematicians who were stuck on the restriction problem connected it to another problem called the Kakeya problem,  a geometry problem about overlapping tubes in space.   The Kakeya problem is not obviously ``half of the restriction problem''-- at first sight,  the Kakeya problem does not even sound like a problem in the same field.   But with work and wisdom,  mathematicians in the field realized that the Kakeya problem is an easier problem that isolates an essential issue in the restriction problem.  

I think many parts of math that sound complicated developed as creative ways to ask more basic questions about places where we were stuck.   
From the outside,  it may look like mathematicians are talking about very complicated things and ``building up'' new theories that  are more complex than the previous theories.   While that isn't wrong,  I think that from the inside it often feels the opposite: we are trying to find more basic questions that lie underneath the previous questions that we can't answer.   Raoul Bott once said,  ``a mathematician is a person who likes to get to the bottom of things''.  



 Coming to a good question is hard work,  full of false starts and crumpled paper just like the rest of the process.  When someone asks a good question,  they have an intuition that trying to answer the question will lead to an interesting idea.   They may not have the interesting idea yet,  but the question is some kind of hint that might lead towards it,  a hint that they can share with other people who can build on it.




\section{Understanding the world}

Mathematics also helps to understand the world,  and especially some of the strange parts of the world.   In my career so far,  I have focused on pure math,  and I feel some imposter syndrome in trying to write about this important issue,  but I will still give it a try.

Mathematics comes partly from trying to understand the world,   partly from looking for connections and relationships,  partly from trying to make ideas more rigorous,  partly from looking for more basic questions,  partly from many other things.    These different aspects somehow help each other.  

When people discovered quantum mechanics,  it was very strange and counterintuitive compared with our everyday experience.    In one sense,  it was extremely different from anything people had seen.   But quantum mechanics can be described mathematically.  And the mathematics involved is actually similar to mathematics that was discovered earlier for very different reasons -- matrices,  inner products, eigenvalues and eigenvectors,  partial differential equations.  

It is true that quantum mechanics is a success story and there are also many less successful stories.  But I still feel that when we encounter something new and strange in the world and try to understand it, math is a tool that can help us.  

One reason I think math can be helpful has to do with asking basic questions.   Trying to study more basic,  more general questions can help to find ideas that may apply in new situations.   When we encounter something new and strange in the world and try to understand it,  it forces us to ask basic questions and rethink basic assumptions.   Sometimes math can help with this because we have already been trying for so long to ask basic questions.  

This also has to do with how things fit together.   There are basic structures that appear in many parts of math.   For instance 
matrices, inner products, eigenvalues and eigenvectors appear in many different parts of math,  including both areas of pure math and areas from all over science and engineering.   

Here is a modest example,  drawn from my own experience as a pure mathematician.   I worked for a long time on an area in Fourier analysis called decoupling,  which addresses questions in pure math about PDE and number theory.   After studying decoupling,  I did some reading about more applied topics in Fourier analysis,  and I noticed many of the same structures.   For instance,  I saw that some calculations pure mathematicians did to help solve a number theory problem were very similar to calculations that scientists did to help recognize helical structures from X-ray diffraction patterns (which helped to discover the double helix structure of DNA).   In another direction,  I saw that the basic structure of decoupling is related to the structure of the fast Fourier transform,  an important algorithm used by many scientists and engineers.

As we keep studying,  things fit together into a larger whole,  which includes math and science and engineering.   There are basic structures that appear in all these areas and connect them.  Cultivating our understanding of these structures helps in the long run to understand the world.

\section{Teaching}

For a long time,  the traditional job of a mathematician has included research and teaching.   It might sound like two jobs smooshed into one,  but that is not how I feel about it.   Everything I've written about is really about research and about teaching,  and I think they are closely related.   

When the teacher and the students discuss math,  things come together (slowly!).    Learning also involves confusion and frustration and doubt and being wrong.   The student and the teacher try to find the confusion.   To make it clearer,  they ask each other questions.   Over time,  they try to find questions that address the most basic issues.    When the math isn't clear,  we say to ourselves,  ``if this part wasn't clear,  there was probably something simpler that wasn't clear either.''     

Sharing ideas this way is a source of joy.    

The teacher explains what they learned to the student.   And the student explains it back in their own words.   And the teacher realizes that their first explanation was a confusing and tries a different way.   And the student sees something that doesn't fit together and rephrases the question.   We keep going,  gradually making it clearer,  gradually seeing how it fits together -- and enjoying the wonder of it all.


The teacher shares what they understand and what they don't understand.   The teacher and the student look out at the strange parts of the world, at things no one understands yet.   Maybe the teacher passes on to the student a question which they think might lead to something interesting.   They wonder what the student will learn someday.

\section{Choosing what to work on}

Each mathematician gets to choose what to work on.

In the beginning part of math education,  students are often told what to do: read these sections, solve these exercises,  do these problems.   But this is not the whole story of learning math.   There is an important moment in a student's math education when they start to make their own time to think about math.   In this personal study time,  they decide what to think about.   In this time,  they start to develop their own interests and their own point of view.

Sometimes, choosing what to work on can be easy.  Sometimes we quickly find a direction we are interested in and we are happy working on it for a long while.  We may not have to worry much about the question of what to work on.

But other times, choosing what to work on can be confusing and frustrating.   We might have to crumple up a lot of paper.   We might have doubts.  This confusion and crumpled paper don't mean that we aren't doing the work right.   The confusions and crumpled paper are part of the work of being a mathematician.

Graduate students and advisors have always worked together on these issues.  Most graduate students, somewhere in their PhDs, go through something hard.   The work stops feeling meaningful.  They may not know why.   The graduate student and advisor talk together and try to figure out what's going on.    Part of the advisor's job is to be there and listen and to help the student find a direction that is meaningful for them. 

And it's worth thinking about these questions even if you aren't feeling frustrated or confused.  

This question of what to work on is not a question that we try to solve so that we can get it out of the way and get on with things.   It's more important than that.    
The poet Rainer Maria Rilke has a famous quote about ``living in the questions''.   The quote has always reminded me of the process of doing math.  Sometimes it reminds me of some mathematical questions,  which a mathematician may spend a long time living with.   But it applies even more to the question of choosing what to work on.

So here is what I suggest.

\begin{itemize}

\item Consider different things to work on and ask yourself,  ``does this feel meaningful to me?''

\item If it does,  do it.

\item If it doesn't,  keep exploring other things.   But don't come to a complete halt while you're waiting to find the perfect direction.   Try a little of this,  and try a little of that.

\item Take time to reflect on what you want to work on.   Brainstorm ideas.   

\item Talk it over with your mentors and with your friends, 

\item Spend some time trying things out and some time reflecting.   Remember that both trying and reflecting are part of the life of a mathematician.

\end{itemize}


\end{document}